\documentclass[10pt,twoside]{amsart}
\usepackage{a4,fancyhdr}
\usepackage{graphicx}
\usepackage{amssymb}
\newtheorem{theorem}{Theorem}[section]
\newtheorem{proposition}{Proposition}[section]
\newtheorem{lemma}{Lemma}[section]

\newtheorem{corollary}{Corollary}[section]

\def\eproof{\hfill $\Box$}

\begin{document}

\title[inverse scattering problem]{\bf Stability of the determination of the surface impedance of an obstacle from the scattering amplitude}

\author{\dag Mourad Bellassoued}
\address{\dag Facult\'e des Sciences de Bizerte, D\'epartement des Math\'ematiques, 7021 Jarzouna Bizerte, Tunisie}
\email{mourad.bellassoued@fsb.rnu.tn}

\author{\ddag Mourad Choulli}
\address{\ddag LMAM, UMR 7122, Universit\'e de Lorraine et CNRS, Ile du Saulcy, 57045 Metz cedex 1, France}
\email{mourad.choulli@univ-lorraine.fr}

\author{\S Aymen Jbalia}
\address{\S LMAM, UMR 7122, Universit\'e de Lorraine et CNRS, Ile du Saulcy, 57045 Metz cedex 1, France\\ et \\ Facult\'e des Sciences de Bizerte, D\'epartement des Math\'ematiques, 7021 Jarzouna Bizerte, Tunisie}
\email{jbalia.aymen@yahoo.fr}

\begin{abstract}
We prove a stability estimate of logarithmic type for the inverse problem consisting in the determination of the surface impedance of an obstacle from the scattering amplitude. We present a simple and direct proof which is essentially based on an elliptic Carleman inequality.
 

\end{abstract}

\maketitle

\tableofcontents

\section{Introduction}
Let $D$ be a bounded subset of $\mathbb{R}^3$. For simplicity, even if it is not always necessary, we assume that $D$ is of class $C^\infty$. As usual, we denote by $u^i$ the incident plan wave : $u^i(x)=e^{ik x\cdot \omega }$, where $k>0$ is the wave number and $\omega \in \mathbb{S}^2$ is the direction of propagation. A simplest model of  the scattering problem for the impedance boundary condition is to find the total field $u=u^i+u^s$, $u^s$ is the scattered field, such that
\begin{eqnarray}
\left\{
\begin{array}{ll}
(\Delta +k^2)u=0 \quad &\mbox{in}\; \mathbb{R}^3\setminus\overline{D},
\\
\partial _\nu u+i\lambda (x)u=0 &\mbox{on}\; \partial D,
\\
\lim\limits_{r\rightarrow \infty}r(\partial_r u^s-iku^s)=0, &r=|x|.
\end{array}
\right.
\end{eqnarray}
Here, $\lambda$ is the surface impedance of the obstacle $D$. The last condition in $(1.1)$ is called the Sommerfeld radiation condition. This condition guarantees that the scattered wave is outgoing.

\begin{theorem}
Assume that $\lambda \in C(D)$ and $\lambda \geq 0$. Then the scattering problem $(1.1)$ has a unique solution $u\in C(\mathbb{R}^3\setminus D)\cap C^2(\mathbb{R}^3\setminus\overline{D})$\footnote{Note that, with this regularity,  Theorem 2.2 in \cite{CK} implies that the solution of $(1.1)$ is analytic in $\mathbb{R}^3\setminus \overline{D}$.}. Moreover, for any $M>0$, there exists a constant $C>0$, depending only on $M$, such that
\begin{equation}
\|u\|_{C(\mathbb{R}^3\setminus D)}\leq C\; \mbox{for any}\; \lambda \in C(D),\; 0\leq \lambda \leq M,
\end{equation}
and for any closed subset $K$ of  $\mathbb{R}^3\setminus \overline{D}$, $\alpha \in \mathbb{N}^3$, there exists a constant $\widetilde{C}$, that can depend only on $K$, $\alpha$ and $M$, such that
\begin{equation}
\|\partial ^\alpha u\|_{C(K)}\leq \widetilde{C}\; \mbox{for any}\; \lambda \in C(D),\; 0\leq \lambda \leq M.
\end{equation}
\end{theorem}

 The proof of the existence and the uniqueness part  is similar to that of Theorem 3.10 in \cite{CK} (Neumann boundary condition) with slight modifications. For sake of completeness, we give a sketch of this proof in Appendix A. In this appendix we give also the proof of estimates $(1.2)$ and $(1.3)$.
 
In order to give a regularity result of the solution of  $(1.1)$, we need to recall the definition of a boundary vector space. Let $m$ be a positive 
integer, $s\in \mathbb{R}$ and $1\leq r\leq \infty$. We consider the vector space
$$
B_{s,r}(\mathbb{R}^m ):=\{ w\in S'(\mathbb{R}^m);\; (1+|\xi |^2)^{s/2}\widehat{w}\in 
L^r(\mathbb{R}^m )\},
$$
where $S'(\mathbb{R}^m)$ is the space of temperate distributions on $\mathbb{R}^m$ and
$\widehat{w}$ is the Fourier transform of $w$. Equipped with the norm 
$$
\|w\|_{B_{s,r}(\mathbb{R}^m )}:=\|(1+|\xi |^2)^{s/2}\widehat{w}\|_{ L^r(\mathbb{R}^m )},
$$
$B_{s,r}(\mathbb{R}^m )$ is a Banach space (it is noted that $B_{s,2}(\mathbb{R}^m )$ is merely the Sobolev
space $H^s(\mathbb{R}^m )$). Using local charts and partition of unity, we construct 
$B_{s,r}(\partial D)$ from $B_{s,r}(\mathbb{R}^{2})$ in the same way as $H^s(\partial D )$ is built 
from $H^s(\mathbb{R}^{2})$.

The space $B_{s,1}(\partial D)$ is very useful because from Theorem 2.1 in \cite{Ch} we know that the multiplication by a function from $B_{s,1}(\partial D)$ defines a bounded operator from $H^s(\partial D)$ into itself.

We proceed similarly to the proof of Theorem 2.3 in \cite{Ch} to prove the $H^\ell$-regularity of the solutions of $(1.1)$. We obtain with the help of the usual elliptic $H^\ell$-regularity (e.g. \cite{LM})  and estimate $(1.3)$ the following  theorem.

\begin{theorem}
Let $\widetilde{\omega} \supset \supset D$ be a bounded $C^\infty$ open subset of $\mathbb{R}^3$, $\ell \geq 0$ an integer and set $\omega =\widetilde{\omega}\setminus D$. If $\lambda \in B_{\ell +1/2,1}(\partial D)\cap C(\partial D)$, then $u$, the solution of the scattering problem $(1.1)$, belongs to $H^{2+\ell}(\omega )$. In addition, for any $M>0$, there exists a constant $C>0$, depending only on $M$ and $\ell$, such that
\begin{equation}
\|u\|_{H^{2+\ell}(\omega )}\leq C\; \mbox{for any}\; \lambda \in B_{\ell +1/2,1}(\partial D)\cap C(\partial D),\; 0\leq \lambda \leq M\; \mbox{and}\; \|\lambda \|_{B_{\ell +1/2,1}(\partial D)}\leq M.
\end{equation}

\end{theorem}
 
\medskip
Since $u^s$ is a radiating solution to the Helmholtz equation, it follows from Theorem 2.5 in \cite{CK} that it has the asymptotic behavior  of an outgoing spherical wave:
$$
u^s(x)=\frac{e^{ik|x|}}{|x|}\left[ u_\infty (\widehat{x})+O\left( \frac{1}{|x|}\right)\right],\quad |x|\rightarrow \infty ,
$$
uniformly in all directions $\widehat{x}=x/|x|$. The function $u_\infty$ defined on $\mathbb{S}^2$ is known as the scattering amplitude or the far field pattern.

In the present paper, we investigate the stability issue of the inverse problem consisting in the determination of the surface impedance coefficient $\lambda$ from the scattering amplitude $u_\infty$. Before stating our main theorem, we need to introduce some geometric assumptions on the domain.

First, assume that $\Omega = \mathbb{R}^3\setminus D$ has the following uniform exterior sphere-interior cone property. 

\medskip
{\bf (GA1)} There exist $\rho >0$ and $\theta \in ]0,\pi /2[$ with the property that for all $\widetilde{x}\in \partial \Omega $, we find $x'\in \mathbb{R}^3\setminus \overline{\Omega}$ such that $B(x',\rho )\cap D =\emptyset$, $B(x',\rho )\cap \overline{D}=\{\widetilde{x}\}$ and
$$
\mathcal{C}(\widetilde{x})=\left\{x\in \mathbb{R}^3;\; (x-\widetilde{x})\cdot \xi >|x-\widetilde{x}|\cos \theta \right\}\subset \Omega ,\; \mbox{where}\; \xi =\frac{\widetilde{x}-x'}{|\widetilde{x}-x'|}.
$$

Set $\mathcal{B}(\widetilde{x},r)=B(x', \rho +r)$ and make the following second geometric assumption on $\Omega$.

\medskip
{\bf (GA2)} There exist $C>0$, $0<\kappa <1$ and $0<r_0$ such that for all $\widetilde{x}\in \partial \Omega$ and $0<r\leq r_0$,
\[
\mathcal{B}(\widetilde{x},r)\cap \partial \Omega \subset B(\widetilde{x},Cr^\kappa )\cap \partial \Omega .
\]

Fix $s>0$, $0<\rho <s$ and let $x'=(0_{\mathbb{R}^2},-\rho )\in \mathbb{R}^3$. A straightforward computation leads
\[
B(x',\rho +r)\cap B((0_{\mathbb{R}^2},-s),s)\subset B\left(0_{\mathbb{R}^3},\left(\frac{s(2\rho +r)}{s-\rho}\right)r\right).
\]
From this simple observation, we deduce that if $D$ has the uniform interior sphere property (or equivalently $\Omega$ has the uniform exterior sphere property) and, under a rigid transform, $\widetilde{x}=0$ and $\partial D$ is represented near $0$ by  the graph $\{x=(y,x_3);\; y\in V,\; x_3=\psi (y)\}$, where $V$ is a neighborhood of $0_{\mathbb{R}^2}$, $\psi$ is such that $\psi (0)=0$ and $\psi \leq 0$\footnote{Note that no regularity assumption is required for $\psi$}, then there exist $C>0$ and $r_0>0$ such that
\[
\mathcal{B}(\widetilde{x},r)\cap \partial \Omega \subset B(\widetilde{x}, Cr )\cap \partial \Omega ,\; 0<r\leq r_0.
\]
Therefore, {\bf (GA2)} is satisfied for instance if $D$ is a $C^2$ convex bounded subset of $\mathbb{R}^3$ (and more generally of $\mathbb{R}^n$). Note that {\bf (GA1)} is automatically satisfied when $D$ is convex.

\begin{theorem} 
Let $M>0$ and $0<\alpha \leq 1$, there exist $C>0$,  $\epsilon >0$, and $\sigma >0$ such that for all $\lambda$, $\widetilde{\lambda}\in B_{3/2,1}(\partial D)\cap C^\alpha (\partial D)$ satisfying\footnote{Here, $\|\cdot \|_{B_{3/2,1}(\partial D)\cap C^\alpha (\partial D)}=\|\cdot \|_{B_{3/2,1}(\partial D)}+\|\cdot \|_{C^\alpha (\partial D)}$.}
$$
\| \lambda \|_{B_{3/2,1}(\partial D)\cap C^\alpha (\partial D)}+\|\widetilde{\lambda}\|_{B_{3/2,1}(\partial D)\cap C^\alpha (\partial D)}\leq M
$$
and
$$
\|\lambda -\widetilde{\lambda}\|_{C(\partial D)} \leq \epsilon ,
$$
we have 
$$
\|\lambda -\widetilde{\lambda}\|_{C(\partial D)}\leq C\left|\ln \left(\frac{\ln |\ln \delta |^2}{|\ln \delta |}\right)\right|^{-\sigma} .
$$
where, $\delta =\|u_\infty (\lambda )-u_\infty (\widetilde{\lambda})\|_{L^2(\mathbb{S}^2)}.$
\end{theorem}

The estimate in the last theorem seems unusual in comparison with the most classical result for inverse elliptic problems which are of log type. The form of the function in the right hand of the last inequality was derived from the usual estimate of the near field by the far field. The reason is to have a simple statement of our stability estimate. However one can rewrite the previous theorem by keeping the original estimate of the near field by far the field (this estimate is given in the proof of  Theorem 1.3 in section 5).

\smallskip
Our result can be seen as an extension of an earlier result by C. Labreuche \cite{La} corresponding to the case where the impedance $\lambda$ is analytic. Similarly to \cite{La}, our proof is based on  a lower bound for the $L^2$-norm of the solution of the scattering problem in any ball around a boundary point. The crucial step in the proof consists in establishing the dependence of the lower bound on the radius of each ball. 

\smallskip
We mention that a result of the same kind as ours was proved by E. Sincich \cite{Si}. The main ingredient in the approach of \cite{Si} is a boundary version of the so-called $A_p$-weight. We observe that $A_p$-weight is also an efficient tool for controlling lower bounds of solutions of elliptic partial differential equations. 

\smallskip
We develop in the present  work a simple  and a direct method which relies essentially on an elliptic Carleman inequality.

\medskip
We make  some geometric assumptions that seem somehow restrictive. We choose to make these assumptions for a better presentation and because the proofs are more simple. We believe that some of these geometric assumptions can be relaxed.

\smallskip
As we said before, the main tool in our method is a Carleman inequality. Precisely, a version with an explicit dependence on data. This is done in section 2. An intermediate result consisting in a quantitative estimate of continuation from Cauchy data is proved section 3. This result is then used in section 4 to derive an appropriate  lower bound for the $L^2$-norm of the solution of the scattering problem in any ball around a boundary point. The results in sections 2 to 4 are given in an arbitrary dimension greater or equal to two. The last section is devoted to the proof of the stability estimate for our inverse scattering problem.

\smallskip
One can see that these results can be adapted to other problems such as the problem of recovering the corrosion coefficient appearing in some usual corrosion detection problems.

\section{An elliptic Carleman inequality}
\setcounter{equation}{0} 

Let $\Omega$ be a Lipschitz bounded domain of $\mathbb{R}^n$ ($n\geq 2$), with boundary $\Gamma$, and let $0\leq \psi \in C^2(\overline{\Omega })$ be  such that 
$$m=\min \left(1,\min_{\overline{\Omega}}|\nabla \psi |\right)>0.$$ 
Fix  $M\geq \max \left(\|\psi\|_{C^2(\overline{\Omega})},1\right)$, where
$$
 \|\psi\|_{C^2(\overline{\Omega})}=\sum_{|\alpha |\leq 2}\|\partial ^\alpha \psi\|_{C(\overline{\Omega})},
$$
and set $\varphi =e^{\varrho \psi}$, $\varrho>0$. 

\begin{proposition}
Let $v\in H^2(\Omega)$. Then
\begin{eqnarray*}
&&\int_\Omega e^{2\tau \varphi}\left(m^4\varrho ^4\tau ^3\varphi ^3v^2+m^2\varrho ^2\tau \varphi|\nabla v|^2\right)dx
\\
&&\leq 8\int_\Omega e^{2\tau \varphi}(\Delta v)^2dx+48\int_\Gamma  e^{2\tau \varphi}\left(M^3\varrho ^3\tau ^3\varphi ^3v^2+M\varrho \tau \varphi|\nabla v|^2\right)d\sigma ,\\ 
&& \hskip 8cm \varrho \geq 6M^3/m^4,\; \tau \geq 88M^6/m^4.
\end{eqnarray*}
\end{proposition}

{\bf Proof.} Set $\Phi =e^{-\tau \varphi}$. Then straightforward computations give
\begin{eqnarray*}
&&\nabla \Phi =-\varrho \tau  \varphi \Phi\nabla\psi
\\
&&\Delta \Phi= \varrho ^2\tau ^2 \varphi ^2\Phi |\nabla \psi |^2-\varrho ^2\tau  \varphi \Phi |\nabla \psi |^2-\varrho \tau  \varphi \Phi \Delta \psi .
\end{eqnarray*}

Let $w\in H^2(\Omega)$. Then from formulae above, we deduce
$$
Lw=[\Phi ^{-1}\Delta \Phi ]w=L_1w+L_2w+cw.
$$
Here,
\begin{align*}
L_1w&=aw+\Delta w,
\\
L_2w&= B\cdot \nabla w+bw,
\end{align*}
with
\begin{align*}
a&=a(x,\varrho ,\tau )= \varrho ^2\tau ^2 \varphi ^2 |\nabla \psi |^2,
\\
B&=B(x,\varrho ,\tau )=-2\varrho \tau \varphi \nabla \psi ,
\\
b&=b(x,\varrho ,\tau )=-2\varrho ^2\tau \varphi|\nabla \psi |^2,
\\
c&=c(x,\varrho ,\tau )=-\varrho \tau \varphi \Delta \psi  +\varrho ^2\tau \varphi |\nabla \psi |^2.
\end{align*}

We have
\begin{equation}
\int_\Omega awB\cdot \nabla wdx=\frac{1}{2}\int_\Omega aB\cdot \nabla w^2dx=-\frac{1}{2}\int_\Omega \mbox{div}(aB) w^2dx+\frac{1}{2}\int_\Gamma aB\cdot \nu w^2d\sigma
\end{equation}
and
\begin{align}
\qquad \int_\Omega \Delta wB\cdot \nabla wdx&=-\int_\Omega \nabla w\cdot \nabla (B\cdot \nabla w)dx+\int_\Gamma B\cdot \nabla w\nabla w\cdot\nu d\sigma
\\
&=-\int_\Omega B'\nabla w\cdot\nabla wdx-\int_\Omega\nabla ^2wB\cdot \nabla wdx+\int_\Gamma B\cdot \nabla w\nabla w\cdot\nu d\sigma.\nonumber
\end{align}
Here, $B'=(\partial _iB_j)$ is the jacobian matrix of $B$ and $\nabla ^2w=(\partial ^2_{ij}w)$ is the hessian matrix of $w$.

But,
$$
\int_\Omega B_i\partial ^2_{ij}w\partial _jwdx=-\int_\Omega \partial _jwB_i\partial ^2_{ij}wdx-\int_\Omega\partial _iB_i(\partial _j w)^2dx+\int_\Gamma B_i(\partial _j w)^2\nu _i d\sigma.
$$
Therefore,
\begin{equation}
\int_\Omega\nabla ^2wB\cdot \nabla wdx=-\frac{1}{2}\int_\Omega \mbox{div}(B)|\nabla w|^2dx+\frac{1}{2}\int_\Gamma |\nabla w|^2B\cdot \nu  d\sigma .
\end{equation}

It follows from $(2.2)$ and $(2.3)$,
\begin{equation}
\int_\Omega \Delta wB\cdot \nabla wdx=\int_\Omega \left[-B'+\frac{1}{2}\mbox{div}(B)I\right]\nabla w\cdot\nabla wdx+\int_\Gamma B\cdot \nabla w\nabla w\cdot\nu d\sigma - \frac{1}{2}\int_\Gamma |\nabla w|^2B\cdot \nu  d\sigma .
\end{equation}

As before, an integration by parts gives
$$
\int_\Omega \Delta wbwdx=-\int_\Omega b|\nabla w|^2dx-\int_\Omega w\nabla b\cdot \nabla wdx+\int_\Gamma bw\nabla w\cdot \nu d\sigma .
$$
Then, using the following inequality
$$
-\int_\Omega w\nabla b\cdot \nabla wdx\geq -\int_\Omega (\varrho ^2\varphi )^{-1}|\nabla b|^2w^2dx-\int_\Omega \varrho ^2\varphi |\nabla w|^2dx ,
$$
we obtain
\begin{equation}
\int_\Omega \Delta wbwdx\geq -\int_\Omega (b+\varrho ^2\varphi )|\nabla w|^2dx-\int_\Omega (\varrho ^2\varphi )^{-1}|\nabla b|^2w^2dx+\int_\Gamma bw\nabla w\cdot \nu d\sigma .
\end{equation}

Now a combination of $(2.1)$, $(2.4)$ and $(2.5)$ leads
\begin{equation}
\int_\Omega L_1wL_2wdx -\int_\Omega c^2w^2dx\geq \int_\Omega fw^2dx+\int_\Omega F\nabla w\cdot \nabla w dx+\int_\Gamma g(w)d\sigma ,
\end{equation}
where,
\begin{align*}
f&=-\frac{1}{2}\mbox{div}(aB)+ab-(\varrho ^2\varphi )^{-1}|\nabla b|^2-c^2,
\\
F&=-B'+\Big(\frac{1}{2}\mbox{div}(B)-b-\varrho ^2\varphi \Big)I,
\\
g(w)&=\frac{1}{2}aw^2B\cdot \nu-\frac{1}{2}|\nabla w|^2B\cdot \nu+B\cdot \nabla w\nabla w \cdot \nu+bw\nabla w \cdot \nu .
\end{align*}

From the following elementary inequality $(s-t)^2\geq s^2/2-t^2$, $s$, $t>0$, we obtain
$$
\|Lw\|_2^2\geq (\|L_1w+L_2w\|_2-\|cw\|_2)^2\geq \frac{1}{2}\|L_1w+L_2w\|_2^2-\|cw\|_2^2\geq \int_\Omega L_1wL_2w dx-\int_\Omega c^2w^2dx.
$$
This and  $(2.6)$ imply
\begin{equation}
\|Lw\|_2^2\geq \int_\Omega fw^2dx+\int_\Omega F\nabla w\cdot \nabla w dx+\int_\Gamma g(w)d\sigma .
\end{equation}

By a straightforward computation, we prove
$$
-\frac{1}{2}\mbox{div}(aB)=\varrho ^3\tau ^3\mbox{div}(\varphi ^3|\nabla \psi |^2\nabla \psi)=\varrho ^3\tau ^3(3\varrho \varphi ^3|\nabla \psi |^4+\varphi ^2\mbox{div}(|\nabla \psi |^2\nabla \psi)).
$$
Therefore,
$$
-\frac{1}{2}\mbox{div}(aB)+ab=\varrho ^3\tau ^3(\varrho \varphi ^3|\nabla \psi |^4+\varphi ^3\mbox{div}(|\nabla \psi |^2\nabla \psi)).
$$
Hence,
\begin{equation}
-\frac{1}{2}\mbox{div}(aB)+ab\geq \varrho ^3\tau ^3\varphi ^3(\varrho m^4-3M^3).
\end{equation}

From now, we assume that $\varrho \geq 1$ and $\tau \geq 1$. Using $-c^2\geq -4\varrho ^4\tau ^2\varphi ^2M^4$, we deduce from $(2.8)$
\begin{equation}
-\frac{1}{2}\mbox{div}(aB)+ab-c^2\geq \varrho ^3\tau ^3\varphi ^3(\varrho m^4-3M^3)-12\varrho ^4\tau ^2\varphi ^2M^4.
\end{equation}

Next, we estimate $|\nabla b|^2$. We have
$$
\nabla b=-2\varrho ^2\tau \nabla (\varphi|\nabla \psi |^2)=-2\varrho ^2\tau (\varrho \varphi |\nabla \psi |^2|\nabla \psi +\varphi \nabla|\nabla \psi |^2).
$$
Consequently,
$$
-|\nabla b|^2\geq -10\varrho ^6\tau ^2\varphi ^2M^6.
$$
This and $(2.9)$ yield
$$
-\frac{1}{2}\mbox{div}(aB)+ab-c^2-(\varrho ^2\varphi )^{-1}|\nabla b|^2\geq \varrho ^3\tau ^3\varphi ^3(\varrho m^4-3M^3)-22\varrho ^4\tau ^2\varphi ^3M^6.
$$
That is,
$$
f\geq \varrho ^3\tau ^3\varphi ^3(\varrho m^4-3M^3)-22\varrho ^4\tau ^2\varphi ^3M^6.
$$
Then,
$$
f\geq \frac{1}{2}\varrho ^4\tau ^3\varphi ^3m^4-22\varrho ^4\tau ^2\varphi ^3M^6=\varrho ^4\tau ^2\varphi ^3\left(\frac{1}{2}\tau m^4-22M^6\right),\quad \varrho \geq 6M^3/m^4.
$$
Hence,
\begin{equation}
f\geq \frac{1}{4}\varrho ^4\tau ^3\varphi ^3m^4,\quad \varrho \geq 6M^3/m^4,\quad \tau \geq 88M^6/m^4.
\end{equation}

We have
\begin{equation}
-B'\xi \cdot \xi = 2\varrho \tau (\varrho \varphi |\nabla \psi \cdot \xi |^2+\varphi \nabla ^2\psi \xi \cdot \xi )\geq -2\varrho \tau \varphi M|\xi |^2,\; \xi \in \mathbb{R}^n.
\end{equation}
On the other hand,
\begin{align*}
 \frac{1}{2}\mbox{div}(B)-b-\varrho ^2\varphi &= -\varrho^2\tau \varphi|\nabla \psi |^2-\varrho \varphi \Delta \psi +2\varrho^2\tau \varphi|\nabla \psi |^2-\varrho ^2\varphi 
\\ 
&=\varrho^2\tau \varphi|\nabla \psi |^2-\varrho \varphi \tau\Delta \psi-\varrho ^2\varphi 
\\ 
&\geq \varrho^2\tau \varphi m^2-\varrho \varphi \tau M-\varrho ^2\varphi,\quad \varrho \geq 4M/m^2.
\end{align*}
A combination of this estimate and $(2.11)$ implies
\begin{equation}
F\xi \cdot \xi \geq \frac{1}{4}\varrho^2\tau \varphi m^2,\quad \varrho \geq 6M/m^2,\quad
\tau \geq 4/m^2,\quad \xi \in \mathbb{R}^n,\; |\xi |=1.
\end{equation}

For $g(w)$, we first note that
$$
|bw\nabla w\cdot \nu |=\sqrt{\varrho |\nabla \psi ||b|}|w|\sqrt{(\varrho|\nabla \psi |)^{-1}|b|}|\nabla w\cdot \nu |\leq \varrho |\nabla \psi ||b|w^2+(\varrho|\nabla \psi |)^{-1}|b||\nabla w|^2.
$$
From this inequality, we easily deduce
\begin{equation}
|g(w)|\leq 2(M^3\varrho ^3\tau ^3 \varphi ^3w^2+M\varrho \tau \varphi |\nabla w|^2).
\end{equation}

Finally, $(2.7)$, $(2.10)$, $(2.12)$ and $(2.13)$ yield
\begin{eqnarray}
&&\int_\Omega \left(m^4\varrho ^4\tau ^3\varphi ^3w^2+m^2\varrho ^2\tau \varphi|\nabla w|^2\right)dx
\\
&&\qquad \leq 4\int_\Omega (Lw)^2dx+8\int_\Gamma \left(M^3\varrho ^3\tau ^3\varphi ^3w^2+M\varrho \tau \varphi|\nabla w|^2\right)d\sigma ,\nonumber
\\
&& \hskip 8cm \varrho \geq 6M^4/m^4,\quad \tau \geq 88M^6/m^4.\nonumber
\end{eqnarray}

Let us now apply this inequality to $w=\Phi ^{-1}v$, $v\in H^2(\Omega )$. We have
\begin{align*}
|\nabla w|^2&=|\Phi ^{-1}\nabla v-\Phi ^{-2}v\nabla \Phi |^2\geq \frac{1}{2}\Phi ^{-2}|\nabla v|^2-\Phi ^{-2}|\Phi ^{-1}\nabla \Phi |^2v^2
\\
&\geq\frac{1}{2}\Phi ^{-2}|\nabla v|^2-\Phi ^{-2}\varrho ^2\tau ^2\varphi ^2M^2v^2.
\\
&\geq\frac{1}{2}\Phi ^{-2}|\nabla v|^2-\Phi ^{-2}\varrho ^2\tau ^2\varphi ^3M^2v^2.
\end{align*}
Hence,
\begin{eqnarray*}
&&2\int_\Omega \left(m^4\varrho ^4\tau ^3\varphi ^3w^2+m^2\varrho ^2\tau \varphi|\nabla w|^2\right)dx \qquad
\\
&& \qquad \qquad\geq \int_\Omega \Phi ^{-2}\left(\left[2m^4\varrho ^4\tau ^3-2\varrho ^2\tau ^2M^2\right]\varphi ^3v^2+m^2\varrho ^2\tau \varphi|\nabla v|^2\right)dx.
\end{eqnarray*}
But,
$$
2m^4\varrho ^4\tau ^3-2\varrho ^2\tau ^2M^2=m^4\varrho ^4\tau ^3+m^4\varrho ^4\tau ^3-2\varrho ^2\tau ^2M^2=m^4\varrho ^4\tau ^3+\varrho ^2\tau ^2(\varrho ^2\tau -M^2).
$$
Therefore,
$$
2m^4\varrho ^4\tau ^3-2\varrho ^2\tau ^2M^2\geq m^4\varrho ^4\tau ^3,\quad \varrho \geq 6M^3/m^4,\quad \tau \geq 88M^6/m^4.
$$

With these inequalities in view, we easily deduce from $(2.14)$,
\begin{eqnarray*}
&&\int_\Omega \Phi ^{-2}\left(m^4\varrho ^4\tau ^3\varphi ^3v^2+m^2\varrho ^2\tau \varphi|\nabla v|^2\right)dx
\\
&&\leq 8\int_\Omega \Phi ^{-2}(\Delta v)^2+48\int_\Gamma \Phi ^{-2} \left(M^3\varrho ^3\tau ^3\varphi ^3v^2+M\varrho \tau \varphi|\nabla v|^2\right)d\sigma ,
\\
&& \hskip 8cm \varrho \geq 6M^3/m^4,\; \tau \geq 88M^6/m^4,
\end{eqnarray*}
which is the desired inequality.
\eproof

\medskip
Let $P$  be a partial differential operator of the form
$$
P=\Delta w+A\cdot \nabla +a,
$$
where $A\in L^\infty (\Omega) ^n$, $a\in L^\infty (\Omega )$.

Fix $\Lambda >0$ satisfying
$$
\Lambda \geq 4\max \left(\|A\|_{L^\infty (\Omega) ^n}^2,\|a\|_{L^\infty (\Omega )}^2\right).
$$
Then a straightforward computation shows
$$
\left(\Delta w\right)^2\leq \left(Pw\right)^2+\Lambda \left(w^2+|\nabla w|^2\right)\; \mbox{in}\; \Omega ,\quad \mbox{for any}\; w\in H^2(\Omega ).
$$

Using this inequality, we obtain as an immediate consequence of the previous proposition the following corollary.

\begin{corollary}
For any $v\in H^2(\Omega )$, we have
\begin{eqnarray*}
&&\int_\Omega e^{2\tau \varphi}\left(m^4\varrho ^4\tau ^3\varphi ^3v^2+m^2\varrho ^2\tau \varphi|\nabla v|^2\right)dx
\\
&&\leq 32\int_\Omega e^{2\tau \varphi}(Pv)^2dx+96\int_\Gamma  e^{2\tau \varphi}\left(M^3\varrho ^3\tau ^3\varphi ^3v^2+M\varrho \tau \varphi|\nabla v|^2\right)d\sigma ,
\end{eqnarray*}
for any $$\varrho \geq 6M^3/m^4,\; \tau \geq \max (88M^6, 16\Lambda )/m^4$$
or
$$\varrho \geq \max (6M^3,16\Lambda )/m^4,\; \tau \geq 88M^6/m^4 .$$
\end{corollary}

We shall need also the following consequence of Proposition 2.1.

\begin{corollary}
Let $\widetilde{\Lambda}$ be given. Then for any $v\in H^2(\Omega )$ satisfying
$$
(\Delta v)^2\leq \widetilde{\Lambda}\left(v^2+|\nabla v |^2\right)\; \mbox{in}\; \Omega ,
$$
we have
$$
\int_\Omega e^{2\tau \varphi}\left(m^4\varrho ^4\tau ^3\varphi ^3v^2+m^2\varrho ^2\tau \varphi|\nabla v|^2\right)dx
\leq 96\int_\Gamma  e^{2\tau \varphi}\left(M^3\varrho ^3\tau ^3\varphi ^3v^2+M\varrho \tau \varphi|\nabla v|^2\right)d\sigma ,
$$
for any $$\varrho \geq 6M^3/m^4,\; \tau \geq \max (88M^6, 16\widetilde{\Lambda} )/m^4$$
or
$$\varrho \geq \max (6M^3,16 \widetilde{\Lambda})/m^4,\; \tau \geq 88M^6/m^4.$$
\end{corollary}

\section{A quantitative estimate of continuation from Cauchy data}
\setcounter{equation}{0}

Let $\Omega$ be as in the previous section. That is a bounded Lipschitz domain of $\mathbb{R}^n$ with boundary $\Gamma$. Let $\Gamma _0$ be a closed subset of $\Gamma$ having nonempty interior. We assume that $\Omega$ has the uniform exterior sphere property at any point of $\Gamma _0$:

\medskip
{\bf (GA0)} there exists $\rho >0$ with the property that, for all $\widetilde{x}\in \Gamma _0$, we find $x_0\in \mathbb{R}^n\setminus \overline{\Omega}$ such that $B(x_0,\rho )\cap \Omega =\emptyset$ and $B(x_0,\rho )\cap \overline{\Omega}=\{\widetilde{x}\}$.

\medskip
We shall use the following notations
$$
\mathcal{B}(\widetilde{x},r_1)=B(x_0, \rho +r_1)\, \quad \mathcal{B}(\widetilde{x},r_1,r_2)=\mathcal{B}(\widetilde{x},r_2)\setminus \overline{\mathcal{B}(\widetilde{x},r_1)},\quad \mathcal{B}=\mathcal{B}(\widetilde{x},d),
$$
where $d=\mbox{diam}(\Gamma _0)$.

Henceforth, $P$ is an operator with bounded coefficients of the form
$$
P=\Delta +A\cdot \nabla +a.
$$
Set
$$
\Lambda = 4\max \left(\|A\|_{L^\infty (\Omega )^n},\|a\|_{L^\infty (\Omega )}\right).
$$

\begin{lemma}
There exist two constants $C>0$ and $0<\gamma <1$ with the property that, for any $0<r\leq d$ and any $u\in H^2(\Omega )$ satisfying $Pu=0$ in $\Omega$,
we have the following estimate
$$
Cr^2\|u\|_{L^2(\mathcal{B}(\widetilde{x},\frac{r}{2})\cap\Omega )}\leq \|u\|_{H^1(\Omega )}^{1-\gamma}\left(\|u\|_{L^2(\mathcal{B}(\widetilde{x},r)\cap \Gamma )}+\||\nabla u|\|_{L^2(\mathcal{B}(\widetilde{x},r)\cap \Gamma )}\right)^\gamma .
$$
\end{lemma}

{\bf Proof.} Pick $\widetilde{x}\in \Gamma _0$. Let $x_0$ be as in {\bf (GA0)} and 
$$\psi (x)=\psi _{\widetilde{x}}(x)=\ln \big((\rho +d)^2/|x-x_0|^2\big). $$ 
Then 
$$
|\nabla \psi (x)|=\frac{2}{|x-x_0|}\geq \frac{2}{\rho+d}=m',\; x\in \overline{\mathcal{B}}.
$$
Set $m=\min (1,m')$ and
$$
M=\max _{\widetilde{y}\in \Gamma _0}\left(1,\sum_{|\alpha |\leq 2}\|\partial ^\alpha \psi _{\widetilde{y}}\|_{C(\overline{\mathcal{B}})}\right).
$$
Let $\chi \in C_c^\infty (\mathcal{B}(\widetilde{x},r))$, $\chi =1$ on $\mathcal{B}(\widetilde{x},\frac{3r}{4})$ and $|\partial ^\alpha \chi |\leq Kr^{-|\alpha |}$, $|\alpha |\leq 2$, where $K$ is a constant independent on $r$.

Let $u\in H^2(\Omega )$ satisfying $Pu=0$ in $\Omega$. We apply Corollary 2.1 to $v=\chi u$. For $\varrho =\max (6M^3,16\Lambda )/m^4$ and $\tau \geq \tau _0=88M^6/m^4$, we obtain
\begin{equation}
C\int_{\mathcal{B}(\widetilde{x},\frac{r}{2})\cap \Omega}e^{2\tau \varphi}u^2dx\leq \int_{\mathcal{B}(\widetilde{x},r)\cap \Omega}e^{2\tau \varphi}(Qu)^2dx+\frac{1}{r^2}\int_{\mathcal{B}(\widetilde{x},r)\cap \Gamma}e^{2\tau \varphi}(u^2+|\nabla u|^2)d\sigma .
\end{equation}
Here and in the sequel, $C$ is a generic constant independent on $r$ and
$$
Qu=2\nabla \chi \cdot \nabla u+\Delta \chi u +A\cdot \nabla \chi u.
$$
Using the properties of $\chi$, we easily prove
$$
\int_{\mathcal{B}(\widetilde{x},r)\cap \Omega}e^{2\tau \varphi}(Qu)^2dx\leq \frac{C}{r^4}\int_{\mathcal{B}(\widetilde{x},\frac{3r}{4},r)\cap\Omega}e^{2\tau \varphi}(u^2+|\nabla u|^2)dx.
$$
Therefore, $(3.1)$ implies
$$
Cr^4\int_{\mathcal{B}(\widetilde{x},\frac{r}{2})\cap\Omega}e^{2\tau \varphi}u^2dx\leq  \int_{\mathcal{B}(\widetilde{x},\frac{3r}{4},r)\cap\Omega}e^{2\tau \varphi}(u^2+|\nabla u|^2)dx
+\int_{\mathcal{B}(\widetilde{x},r)\cap \Gamma}e^{2\tau \varphi}(u^2+|\nabla u|^2)d\sigma .
$$

We have
$$
\varphi =e^{\varrho \ln \big((\rho +d)^2/|x-x_0|^2\big)}=\frac{(\rho +d)^{2\varrho}}{|x-x_0|^{2\varrho}}.
$$
Consequently,
\begin{eqnarray}
&& Cr^4e^{2\tau \varphi _0}\int_{\mathcal{B}(\widetilde{x},\frac{r}{2})\cap\Omega}u^2dx\leq e^{2\tau \varphi _1}\int_{\mathcal{B}(\widetilde{x},\frac{3r}{4},r)\cap\Omega}(u^2+|\nabla u|^2)dx
\\
 &&\hskip 6cm +e^{2\tau \varphi _2}\int_{\mathcal{B}(\widetilde{x},r)\cap \Gamma}(u^2+|\nabla u|^2)d\sigma ,\nonumber
\end{eqnarray}
where,
$$
\varphi _0=\frac{(\rho +d)^{2\varrho}}{(\rho +\frac{r}{2})^{2\varrho}},\quad \varphi _1=\frac{(\rho +d)^{2\varrho}}{(\rho +\frac{3r}{4})^{2\varrho}},\quad \varphi _2=\frac{(\rho +d)^{2\varrho}}{\rho ^{2\varrho}}.
$$
By the mean value theorem, for some $\theta \in ]0,1[$,
\begin{align*}
\varphi _0-\varphi _1&= \Big( (\rho +\frac{r}{2})-(\rho +\frac{3r}{4})\Big)\frac{-2\varrho (\rho +d)^{2\varrho}}{(\theta (\rho +\frac{r}{2})+(1-\theta)(\rho +\frac{3r}{4}))^{2\varrho +1}}
\\
&\geq \frac{\varrho r}{2}\frac{(\rho +d)^{2\varrho}}{(\rho +\frac{3r}{4}))^{2\varrho +1}}
\\
&\geq \frac{\varrho r}{2}\frac{(\rho +d)^{2\varrho}}{(\rho +\frac{3d}{4}))^{2\varrho +1}}=\alpha r 
\end{align*}
with
$$
\alpha= \frac{\varrho (\rho +d)^{2\varrho}}{2(\rho +\frac{3d}{4}))^{2\varrho +1}}.
$$
Similarly, we prove
$$
\varphi _2-\varphi _1\leq \beta r, 
$$
with
$$
\beta =\frac{\varrho (\rho +d)^{2\varrho}}{\rho ^{2\varrho +1}}.
$$

We obtain from $(3.2)$, 
$$
Cr^4\int_{\mathcal{B}(\widetilde{x},\frac{r}{2})\cap\Omega}u^2dx\leq e^{-\alpha r\tau }\int_{\mathcal{B}(\widetilde{x},\frac{3r}{4},r)\cap\Omega }(u^2+|\nabla u|^2)dx
 +e^{\beta r\tau }\int_{\mathcal{B}(\widetilde{x},r)\cap \Gamma}(u^2+|\nabla u|^2)d\sigma .
$$
In particular,
\begin{equation}
Cr^4\int_{\mathcal{B}(\widetilde{x},\frac{r}{2})\cap\Omega}u^2dx\leq e^{-\alpha r\tau }\int_\Omega (u^2+|\nabla u|^2)dx
 +e^{\beta r\tau }\int_{\mathcal{B}(\widetilde{x},r)\cap \Gamma}(u^2+|\nabla u|^2)d\sigma .
\end{equation}

Let us introduce the following temporary notations
\begin{align*}
A&=\int_\Omega (u^2+|\nabla u|^2)dx,
\\
I&=\int_{\mathcal{B}(\widetilde{x},r)\cap \Gamma}(u^2+|\nabla u|^2)d\sigma ,
\\
J&=Cr^4\int_{\mathcal{B}(\widetilde{x},\frac{r}{2})\cap\Omega}u^2dx.
\end{align*}
Then, $(3.3)$ becomes
\begin{equation}
J\leq e^{-\alpha r\tau }A+e^{\beta r\tau }I.
\end{equation}

Let
$$
\tau _1=\frac{\ln (A/I)}{\alpha r+\beta r}.
$$
If $\tau _1\geq \tau _0$, then $\tau =\tau _1$ in $(3.4)$ yields
\begin{equation}
J\leq A^{\frac{\alpha r}{\alpha r+\beta r}}I^{\frac{\beta r}{\alpha r+\beta r}}=A^{\frac{\alpha }{\alpha +\beta }}I^{\frac{\beta }{\alpha +\beta }}.
\end{equation}
If $\tau _1<\tau _0$, then $A <e^{\tau _0(\alpha +\beta )r}I\leq e^{\tau _0(\alpha +\beta )d}I$. Since
$$
J=Cr^4\int_{\mathcal{B}(\widetilde{x},\frac{r}{2})\cap\Omega}u^2dx\leq Cd^2A,
$$
we have,
\begin{equation}
J\leq CI=CI^{\frac{\alpha }{\alpha +\beta }}I^{\frac{\beta }{\alpha +\beta }}\leq CA^{\frac{\alpha }{\alpha +\beta }}I^{\frac{\beta }{\alpha +\beta }}.
\end{equation}
Hence, in any case, one of estimates $(3.5)$ and $(3.6)$ holds. That is, in terms of our original notations,
$$
Cr^2\|u\|_{L^2(\mathcal{B}(\widetilde{x},\frac{r}{2})\cap\Omega )}\leq \|u\|_{H^1(\Omega )}^{1-\gamma}\left(\|u\|_{L^2(\mathcal{B}(\widetilde{x},r)\cap \Gamma )}+\| |\nabla u|\|_{L^2(\mathcal{B}(\widetilde{x},r)\cap \Gamma )}\right)^\gamma ,
$$
with $\gamma =\frac{\beta}{\alpha + \beta}$.
The proof is then complete.
\eproof

\begin{corollary}
There exist two constants $C>0$ and $0<\gamma <1$ with the property that, for any $0<r\leq d$ and any $u\in H^2(\Omega )$ satisfying $Pu=0$, we have the following estimates
\begin{equation}
Cr^2\|\nabla u\|_{L^2(\mathcal{B}(\widetilde{x},\frac{r}{4})\cap\Omega )}\leq \|u\|_{H^2(\Omega )}^{1-\gamma /2}\left(\|u\|_{L^2(\mathcal{B}(\widetilde{x},r)\cap \Gamma )}+\||\nabla u|\|_{L^2(\mathcal{B}(\widetilde{x},r)\cap \Gamma )}\right)^{\gamma /2} 
\end{equation}
and
\begin{equation}
Cr^2\|u\|_{H^1(\mathcal{B}(\widetilde{x},\frac{r}{4})\cap\Omega )}\leq \|u\|_{H^2(\Omega )}^{1-\gamma/2}\left(\|u\|_{L^2(\mathcal{B}(\widetilde{x},r)\cap \Gamma )}+\||\nabla u|\|_{L^2(\mathcal{B}(\widetilde{x},r)\cap \Gamma )}\right)^{\gamma /2}.
\end{equation}
\end{corollary}

{\bf Proof.} Pick $\chi \in C_c^\infty (\mathcal{B}(\widetilde{x},\frac{r}{2}))$ satisfying $\chi =1$ in $\overline{\mathcal{B}(\widetilde{x},\frac{r}{4})}$ and $|\partial ^\alpha \chi |\leq Kr^{-|\alpha |}$, $|\alpha |\leq 2$, where $K$ is a constant indepedent on $r$. Let $u\in H^2(\Omega )$ satisfying $Pu=0$. From the usual interpolation inequalities, there exists a constant $c=c(\Omega )>0$ such that
$$
\||\nabla (\chi u)|\|_{L^2(\Omega )}\leq c\|\chi u\|_{L^2(\Omega )}^{1/2}\|\chi u\|_{H^2(\Omega )}^{1/2}.
$$
Hence,
\begin{equation}
\||\nabla u|\|_{L^2(\mathcal{B}(\widetilde{x},\frac{r}{4})\cap\Omega )}\leq cr^{-1}\|u\|_{L^2(\mathcal{B}(\widetilde{x},\frac{r}{2})\cap\Omega )}^{1/2}\|u\|_{H^2(\Omega )}^{1/2}.
\end{equation}
On the other hand, it follows from Lemma 3.1
\begin{equation}
Cr^2\|u\|_{L^2(\mathcal{B}(\widetilde{x},\frac{r}{2})\cap\Omega )}\leq \|u\|_{H^2(\Omega )}^{1-\gamma}\left(\|u\|_{L^2(\mathcal{B}(\widetilde{x},r)\cap \Gamma )}+\|\nabla u\|_{L^2(\mathcal{B}(\widetilde{x},r)\cap \Gamma )}\right)^\gamma  .
\end{equation}
Therefore, $(3.7)$ is a consequence of $(3.9)$ and $(3.10)$.

Next, as the trace mapping
$$
w\in H^2\rightarrow (w,\nabla w)\in L^2(\Gamma )^{n+1}
$$
is bounded and
$$
\|u\|_{L^2(\mathcal{B}(\widetilde{x},r)\cap\Gamma )}+\||\nabla u|\|_{L^2(\mathcal{B}(\widetilde{x},r)\cap\Gamma )}\leq \|u\|_{L^2(\Gamma )}+\||\nabla u|\|_{L^2(\Gamma )},
$$
we have
$$
\|u\|_{L^2(\mathcal{B}(\widetilde{x},r)\cap\Gamma )}+\||\nabla u|\|_{L^2(\mathcal{B}(\widetilde{x},r)\cap\Gamma )}\leq K'\|u\|_{H^2(\Omega)}.
$$
Here $K'$ is a constant independent on $r$.

This estimate in $(3.10)$ yields
\begin{equation}
Cr^2\|u\|_{L^2(\mathcal{B}(\widetilde{x},\frac{r}{2})\cap \Omega )}\leq \|u\|_{H^2(\Omega )}^{1-\gamma /2}\left(\|u\|_{L^2(\mathcal{B}(\widetilde{x},r)\cap \Gamma )}+\||\nabla u|\|_{L^2(\mathcal{B}(\widetilde{x},r)\cap \Gamma )}\right)^{\gamma /2}.
\end{equation}
We complete the proof by noting that $(3.8)$ follows from a combination of $(3.7)$ and $(3.11)$.
\eproof

\section{Lower bound for solutions of elliptic equations}
\setcounter{equation}{0}

As in the previous section, $\Omega$ is a bounded Lipschitz domain of $\mathbb{R}^n$, with boundary $\Gamma$, and $P$ is an elliptic operator of the form
$$
P=\Delta +A\cdot \nabla +a,
$$
with $A\in L^\infty(\Omega )^n$ and $a\in L^\infty (\Omega)$. Let
$$
\Lambda =4\max \big(\|A\|^2_{L^\infty (\Omega )^n}+\|a\|^2_{L^\infty (\Omega )} \big).
$$

\medskip
We start with a three sphere inequality. Set $B(i)=B(0,i)$, $i=1,2,3$ and $r_0=\frac{1}{3}\mbox{diam}(\Omega )$. Let $y\in \Omega$ and $u\in H^1(B(y,3r))$, where $0<r< \frac{1}{3}\mbox{dist}(y,\Gamma ) (\leq r_0)$. If
$$
v(x)=u(rx+y),\; x\in B(3),
$$
a simple change of variables leads to the following inequalities
\begin{equation}
c_\ast r^{1-n/2}\|u\|_{H^1(B(y,ir))}\leq \|v\|_{H^1(B(i))}\leq c^\ast r^{-n/2}\|u\|_{H^1(B(y,ir))}.
\end{equation}
Here,
$$
c_\ast =\min (1,r_0), \quad c^\ast =\max (1,r_0).
$$
In addition, if  $u$ satisfies $Pu=0$ in $B(y,3r)$, then a straightforward computation yields
$$
(\Delta v)^2\leq \widetilde{\Lambda} (v^2+|\nabla v|^2)\; \mbox{in}\; B(3),
$$
where, $\widetilde{\Lambda}=\Lambda r_0^2\max (1,r_0^2)$.

We apply Corollary 2.2 to $w=\chi v$, where $\chi \in C_c^\infty (U)$, $\chi =1$ in $K$, with
$$
U=\{x\in \mathbb{R}^n;\; 1/2<|x|<3\},\quad K=\{x\in \mathbb{R}^n;\; 1\leq r\leq 5/2\}.
$$
Similarly to the previous section, we prove the following three spheres inequality in which the constant $C>0$ and $0<\alpha <1$ depend only on $\Lambda$ and $r_0$.
\begin{equation}
\|v\|_{H^1(B(2))}\leq C\|v\|^\alpha _{H^1(B(1))}\|v\|^{1-\alpha} _{H^1(B(3))}.
\end{equation}

The following lemma is a consequence of $(4.1)$ and $(4.2)$.

\begin{lemma}
There exist $C>0$ and $0<\alpha <1$, depending only on $\Lambda >0$ and $r_0>0$ such that, for all  $u\in H^2(\Omega )$ satisfying $Pu=0$ in $\Omega$, $y\in \Omega$ and $0<r\leq \frac{1}{3 }\mbox{dist}(y,\Gamma )$,  
$$
r\|u\|_{H^1(B(y,2r))}\leq C\|u\|_{H^1(B(y,r))}^\alpha\|u\|_{H^1(B(y,3r))}^{1-\alpha}.
$$
\end{lemma}

Next, let $\Omega _0=\mathbb{R}^n\setminus K$, where $K$ is a compact subset of $\mathbb{R}^n$ with nonempty interior. For simplicity, we assume that $0\in \mathring{K}$. We also make the assumption that $\Omega _0$ has the uniform exterior sphere-interior cone property {\bf (GA1)}. Recall that this geometric assumption is the following one.

\medskip
{\bf (GA1)} There exist $\rho >0$ and $\theta \in ]0,\pi /2[$ with the property that, for all $\widetilde{x}\in \partial \Omega _0$, we find $x'\in \mathbb{R}^n\setminus \overline{\Omega} _0$ such that $B(x',\rho )\cap \Omega _0=\emptyset$, $B(x',\rho )\cap \overline{\Omega}_0=\{\widetilde{x}\}$ and
$$
\mathcal{C}(\widetilde{x})=\{x\in \mathbb{R}^n;\; (x-\widetilde{x})\cdot \xi >|x-\widetilde{x}|\cos \theta \}\subset \Omega _0,\; \mbox{where}\; \xi =\frac{\widetilde{x}-x'}{|\widetilde{x}-x'|}.
$$

Fix $R>4\sup_{K}|x|$ and set $\Omega =\Omega _0\cap B(0,R)$.

\begin{theorem}
Let $M>0$, $\tau >0$ and $\Lambda _0>0$ be given. There exist $C>0$, $\eta >0$ and $r^\ast >0$ such that for all $u\in H^{5/2}(\Omega )$ satisfying
\begin{eqnarray*}
\left\{
\begin{array}{lll}
Pu=0\;\; \mbox{in}\;\; \Omega ,
\\
|\partial _\nu u|\leq \Lambda _0|u| \;\; \mbox{on}\;\; \Gamma _0,
\\
|u|\geq \tau \; \mbox{in}\;\; \Omega \cap \{x\in \mathbb{R}^n;\; |x|\geq R/12\} ,
\\
\|u\|_{H^{5/2}(\Omega )}\leq M,
\end{array}
\right.
\end{eqnarray*}
$\widetilde{x}\in \Gamma _0=\partial \Omega _0$ and $0<r\leq \mbox{diam}(\Gamma _0)$, we have,
$$
e^{-\frac{C}{r^\eta}}\leq \|u\|_{L^2(\mathcal{B}(\widetilde{x},r)\cap \Gamma _0)}.
$$
\end{theorem}

{\bf Proof.} Let $\widetilde{x}\in \Gamma _0$ and $x'\in \mathbb{R}^n\setminus \overline{\Omega}$ such that $B(x',\rho )\cap \Omega =\emptyset$, $B(x',\rho )\cap \overline{\Omega}=\{\widetilde{x}\}$. Set, for $0<r<\mbox{diam} (\Gamma _0)$,
$$
x''=\widetilde{x}+r\xi ,\quad x_0=\frac{\widetilde{x}+x''}{2}.
$$
Clearly, $B(x_0,r/2)\subset \mathcal{B}(\widetilde{x},r)\cap\Omega$ and, for $d_0=|x_0-\widetilde{x}|$, $\rho _0=(d_0\sin \theta )/3$, $B(x_0,3\rho _0)\subset \mathcal{C}(\widetilde{x})$.

By induction in $k$, we construct a sequence of balls $(B(x_k, 3\rho _k))$, contained in $\mathcal{C}(\widetilde{x})$, as follows
\begin{eqnarray*}
\left\{
\begin{array}{ll}
x_{k+1}=x_k+\alpha _k \xi ,
\\
\rho_{k+1}=\mu \rho_k ,
\\
d_{k+1}=\mu d_k,
\end{array}
\right.
\end{eqnarray*}
where,
$$
d_k=|x_k-\widetilde{x}|,\;\; \rho _k=cd_k,\;\; \alpha _k=(1-\mu)d_k ,
$$
with
$$ 
c=\frac{\sin \theta}{3},\;\; \mu =\frac{3+2\sin \theta}{3+\sin \theta}.
$$
Note that our construction of these balls guaranties  that
\begin{equation}
B(x_{k+1},\rho _{k+1})\subset B(x_k,2\rho _k).
\end{equation}

Let us denote by $N$ the smallest integer such that $d_N\geq R/8$. Since $d_N=\mu ^N\frac{r}{2}$, 
\begin{equation}
\frac{\ln \frac{R}{4r}}{\ln \mu}\leq N <\frac{\ln \frac{R}{4r}}{\ln \mu}+1
\end{equation}
or equivalently,
$$
N=\left[ \frac{\ln \frac{R}{4r}}{\ln \mu}\right].
$$
If $0\leq k \leq N$, then
$$
|x_k|+3\rho _k\leq |\widetilde{x}|+d_N+\sin \theta d_N\leq R/4+\mu R/4 \leq 3R/4.
$$
Here,  we used that $1<\mu <2$ and $d_N=\mu \mu^{N-1}d_0<\mu \frac{R}{8}$. 

Also, for $x\in B(x_N,\rho_N)$,
$$
|x|\geq |x_N|-\rho _N\geq  d_N-\frac{\sin \theta}{3}d_N\geq \frac{2}{3}d_N\geq \frac{R}{12}.
$$

In other words,
\begin{equation}
B(x_k,3\rho _k)\subset \Omega ,\; 0\leq k \leq N\;  \mbox{and}\; B(x_N,\rho _N)\subset \Omega \cap \{x\in \mathbb{R}^n;\; |x|\geq R/12\}.
\end{equation}

We obtain by applying Lemma 4.1,
$$
\rho_ 0\|u\|_{H^1(B(x_0,2\rho _0))}\leq CM^{1-\alpha}\|u\|_{H^1(B(x_0,\rho _0))}^\alpha .
$$
But from $(4.3)$, $B(x_1,\rho _1)\subset B(x_0,2\rho _0)$. Therefore
\begin{equation}
\rho_ 0\|u\|_{H^1(B(x_1, \rho _1))}\leq CM^{1-\alpha}\|u\|_{H^1(B(x_0,\rho _0))}^\alpha .
\end{equation}

Set
$$
I_k=\|u\|_{H^1(B(x_k,\rho _k))}.
$$
Then $(4.6)$ can be rewritten as follows
$$
I_1\leq \frac{C}{\rho _0}M^{1-\alpha}I_0^\alpha .
$$

Using an induction in $k$, we prove
$$
I_k\leq \frac{C^{1+\alpha +\ldots +\alpha ^{k-1}}}{\rho _{k-1}\rho_{k-2}^\alpha \ldots \rho_0^{\alpha ^{k-1}}}I_0^{\alpha ^k}M^{(1-\alpha )(1+\alpha +\ldots +\alpha ^{k-1})}.
$$
From the inequality 
$$
\rho _{k-1}\rho_{k-2}^\alpha \ldots \rho_0^{\alpha ^{k-1}}=\mu ^m\rho_0^{\frac{1-\alpha ^k}{1-\alpha}}\geq \rho_0^{\frac{1-\alpha ^k}{1-\alpha}},\; \mbox{with}\; m=\sum_{j=0}^{k-2}(k-1-j)\alpha ^j,
$$
it follows
\begin{equation}
I_k\leq \left(\frac{C}{\rho _0}\right)^{\frac{1-\alpha ^k}{1-\alpha}}M^{1-\alpha ^k}I_0^{\alpha ^k}.
\end{equation}

We have
$$
\frac{C}{\rho _0}=\frac{6C}{(\sin \theta )r}.
$$
Hence, we find $r^\ast >0$ such that
$$
\frac{6C}{(\sin \theta )r}\geq 1,\; \mbox{if}\; 0<r\leq r^\ast.
$$

From now we assume that $0<r\leq r^\ast$. We derive from $(4.7)$,
\begin{equation}
I_N \leq\widetilde{M} \left(\frac{C}{r}\right)^\beta I_0^{\alpha ^N},
\end{equation}
where, 
$$
\beta = \frac{1}{1-\alpha},\quad \widetilde{M}=\max (1,M).
$$

Now as $|u|\geq \tau$ in $B(x_N,\rho _N)$, we have
$$
I_N\geq \tau \left|\mathbb{S}^{n-1}\right|^{1/2}\rho _N^{n/2}=\tau \left|\mathbb{S}^{n-1}\right|^{1/2}(\mu ^N\rho _0)^{n/2}
$$
and since $\mu >1$, we deduce,
\begin{equation}
I_N\geq Cr^{n/2} .
\end{equation}

A combination of $(4.8)$ and $(4.9)$ leads
$$
Cr^{\gamma}\leq I_0^{\alpha ^N},\; \mbox{with}\; \gamma=n/2+\beta.
$$
That is,
\begin{equation}
(Cr)^{\gamma /\alpha ^N}\leq I_0.
\end{equation}

By $(4.4)$, we have
$$
\frac{1}{\alpha ^N}=e^{N|\ln \alpha|}<e^{|\ln \alpha|(\ln R+4|\ln r|+1)}.
$$
Therefore, reducing $r^\ast$ if necessary, 
$$
\frac{1}{\alpha ^N}=e^{N|\ln \alpha|}<e^{6|\ln \alpha||\ln r|}=\frac{1}{r^s},\; \mbox{with}\; s=6|\ln \alpha|.
$$

Reducing once again $r^\ast$ if necessary, we  assume that $Cr<1$ in $(4.10)$ (for any $0<r\leq r^\ast$). Then,
$$
(Cr)^{\gamma /\alpha ^N}\geq (Cr)^{\gamma /r^s}=e^{-\frac{\gamma}{r^s}\ln (\frac{1}{Cr})}\geq r^{-\frac{\gamma}{Cr^{s+1}}}.
$$
This and $(4.10)$ imply, where $\eta =s+1$,
$$
e^{-\frac{C}{r^\eta}}\leq \|u\|_{H^1(B(x_0,\rho _0))}\leq \|u\|_{H^1(\mathcal{B}(\widetilde{x},r)\cap\Omega )}.
$$
Combined with $(3.8)$, this estimate yields
\begin{equation}
e^{-\frac{C}{r^\eta}}\leq \|u\|_{L^2(\mathcal{B}(\widetilde{x},r)\cap \Gamma _0)}+\|\nabla u\|_{L^2(\mathcal{B}(\widetilde{x},r)\cap \Gamma _0)}.
\end{equation}
According to our assumption,
$$
|\nabla u|^2= (\partial _\nu u)^2+|\nabla _\tau u|^2\leq \max \left(1,\Lambda )(u^2+|\nabla _\tau u|^2\right).
$$
Hence, $(4.11)$ implies
\begin{equation}
e^{-\frac{C}{r^\eta}}\leq \|u\|_{H^1(\mathcal{B}(\widetilde{x},r)\cap \Gamma _0)}.
\end{equation}

We now estimate the $H^1$ norm in the right hand of the previous inequality by an $L^2$ norm with the help of an interpolation inequality. To this end, let $\chi \in C_c^\infty (\mathcal{B}(\widetilde{x},2r))$ satisfying $\chi =1$ in $\mathcal{B}(\widetilde{x},r)$ and $|\partial ^\alpha \chi |\leq Cr^{-|\alpha |}$, $|\alpha |\leq 2$.

From classical interpolation inequalities, it follows
\begin{eqnarray*}
 \|u\|_{H^1(\mathcal{B}(\widetilde{x},r)\cap \Gamma _0)}&\leq&  \|\chi u\|_{H^1(\Gamma _0)}\leq C\|\chi u\|_{H^2(\Gamma _0)}^{1/2}\|\chi u\|_{L^2(\Gamma _0)}^{1/2}
\\
&\leq & Cr^{-1}\|u\|_{H^2(\Gamma _0)}^{1/2}\|u\|_{L^2(\mathcal{B}(\widetilde{x},2r)\cap\Gamma _0)}^{1/2}.
\end{eqnarray*}
On the other hand, by classical trace theorems, we have
$$
\|u\|_{H^2(\Gamma _0)}\leq \|u\|_{H^{5/2}(\Omega )}\leq M.
$$
The last two estimates together with $(4.12)$ lead to the desired inequality.
\eproof

\medskip
It the sequel, we assume that $\Omega _0$ possesses ${\bf (GA1)}$ and

\medskip
{\bf (GA2)} there exist $C>0$, $0<\kappa <1$ and $0<r_0$ such that, for all $\widetilde{x}\in \Gamma _0$ and $0<r\leq r_0$,
$$
\mathcal{B}(\widetilde{x},r)\cap \Gamma _0\subset B(\widetilde{x},Cr^\kappa )\cap \Gamma _0.
$$

\medskip
Under this new geometric assumption, we deduce from  Theorem 3.1 the following corollary.

\begin{corollary}
Let $M>0$, $\tau >0$ and $\Lambda _0 >0$ be given. There exist $C>0$, $\eta >0$ and $r^\ast >0$ such that, for all $\widetilde{x}\in \Gamma _0=\partial \Omega _0$, $0<r\leq r_0$ and $u\in H^{5/2}(\Omega )$ satisfying
\begin{eqnarray*}
\left\{
\begin{array}{lll}
Pu=0\;\; \mbox{in}\;\; \Omega ,
\\
|\partial _\nu u|\leq \Lambda _0|u|\;\; \mbox{on}\;\; \Gamma _0 ,
\\
|u|\geq \tau \;\; \mbox{in}\;\; \Omega \cap \{x\in \mathbb{R}^n;\; |x|\geq R/12\} ,
\\
\|u\|_{H^{5/2}(\Omega )}\leq M,
\end{array}
\right.
\end{eqnarray*}
we have,
$$
e^{-\frac{C}{r^\eta}}\leq \|u\|_{L^2(B(\widetilde{x},r)\cap \Gamma _0)}.
$$
\end{corollary}

Next, we derive a result on which is based our stability estimate for the inverse problem consisting in the determination of the surface impedance of an obstacle in terms of boundary Cauchy data. Recall that $f\in C^\alpha (\Gamma _0)$ if there exists $L\geq  0$ such that
\begin{equation}
|f(x)-f(x')|\leq L |x-y|^\alpha ,\; x,\, x'\in \Gamma _0.
\end{equation}
We denote by $[f]_\alpha$ the infimum of $L$'s for which  $(4.13)$ is satisfied.

\begin{proposition}
Let $M>0$, $\tau >0$, $0<\alpha \leq 1$ and $\Lambda  _0>0$ be given. There exist $C>0$, $\epsilon >0$ and $\sigma >0$ such that, for all  $u\in H^{5/2}(\Omega )$ satisfying
\begin{eqnarray}
\left\{
\begin{array}{lll}
Pu=0\;\; \mbox{in}\; \Omega ,
\\
|\partial _\nu u|\leq \Lambda _0|u|\;\; \mbox{on}\;\; \Gamma _0 ,
\\
|u|\geq \tau \;\; \mbox{in}\;\; \Omega \cap \{x\in \mathbb{R}^n;\; |x|\geq R/12\} ,
\\
\|u\|_{H^{5/2}(\Omega )}\leq M
\end{array}
\right.
\end{eqnarray}
and, for all $f\in C^\alpha (\Gamma _0)$ satisfying $[f]_\alpha\leq M$, $\|f\|_{L^\infty (\Gamma _0)}\leq \epsilon$, 
$$
\|f\|_{L^\infty (\Gamma _0)}\leq \frac{C}{\left|\ln \left[\|fu\|_{L^\infty (\Gamma _0)}\right]\right|^\sigma }.
$$
\end{proposition}

The following lemma will be used in the proof of Proposition 4.1. Hereafter, $r^\ast$ is the same as in Corollary 4.1.

\begin{lemma}
There exist $\delta ^\ast$ such that, for all $u$ as in Corollary 4.1 satisfying $u\in C(\Gamma _0)$, $\widetilde{x}\in \Gamma _0$ and $0< \delta \leq \delta ^\ast$,
$$
\{x\in B(\widetilde{x}, r^\ast )\cap \Gamma _0;\; |u(x)|\geq \delta\}\neq\emptyset .
$$
\end{lemma}

{\bf Proof.} Otherwise, we  find a sequence $(\delta _k)$, $0<\delta _k\leq \frac{1}{k}$, $(u_k)$  satisfying the assumptions of Corollary 4.1 with $u_k\in C(\Gamma _0)$, for each $k$,  and $(\widetilde{x}_k )$ in  $\Gamma _0$ such that,
$$
\{x\in B(\widetilde{x}_k, r^\ast )\cap \Gamma _0;\; |u_k(x)|\geq \delta _k\}=\emptyset .
$$
In particular,
$$
|u_k|\leq \frac{1}{k}\; \mbox{in}\; B(\widetilde{x}_k, r^\ast )\cap \Gamma _0.
$$
Therefore, we have, by applying Corollary 4.1,
$$
e^{-\frac{C}{(r^\ast )^\eta}}\leq \frac{1}{k}|B(\widetilde{x}_k, r^\ast )\cap \Gamma _0|
\leq  \frac{1}{k}|\Gamma _0|,\; \mbox{for all}\; k\geq 1,
$$
which is impossible. This leads to the desired contradiction and proves the lemma.
\eproof

\medskip
{\bf Proof of Proposition 4.1.} Let $\delta ^\ast$ be as in the previous lemma, $0<\delta <\delta ^\ast$, $u\in H^{5/2}(\Omega )$ satisfying $(3.14)$  and $f\in C^\alpha (\Gamma _0)$.

Let $\widetilde{x}\in \Gamma _0$. If $|u(\widetilde{x} )|\geq \delta$ then
\begin{equation}
|f(\widetilde{x})|\leq \frac{1}{\delta}|f(\widetilde{x})u(\widetilde{x})|.
\end{equation}

Let $\widetilde{x}\in \Gamma _0$ such that $|u(\widetilde{x})|<\delta$ and set
$$
r=\sup \{0<\rho ;\; |u|<\delta \; \mbox{on}\; B(\widetilde{x},\rho )\cap \Gamma _0\}.
$$
From Lemma 4.2 , we know that
$$
\{x\in B(\widetilde{x}, r^\ast )\cap \Gamma _0;\; |u(x)|\geq \delta\}\neq\emptyset .
$$
Hence, $r \leq r^\ast$ and 
$$
\partial B(\widetilde{x},r)\cap \{x\in B(\widetilde{x}, r^\ast )\cap \Gamma _0;\; |u(x)|\geq \delta\}\neq\emptyset .
$$

Let $\widehat{x}\in\partial B(\widetilde{x},r)$ be such that  $|u(\widehat{x})|\geq \delta $. We have,
$$
|f(\widetilde{x})|\leq |f(\widetilde{x})-f(\widehat{x})|+|f(\widehat{x})|\leq [f]_\alpha|\widetilde{x}-\widehat{x}|^\alpha +\frac{1}{\delta}|f(\widehat{x})u(\widehat{x})|
$$
and then,
$$
|f(\widetilde{x})|\leq |f(\widetilde{x})-f(\widehat{x})|+|f(\widehat{x})|\leq Mr^\alpha +\frac{1}{\delta}|f(\widehat{x})u(\widehat{x})|.
$$

This and $(4.15)$ show
\begin{equation}
\|f\|_{L^\infty (\Gamma _0)}\leq Mr^\alpha +\frac{1}{\delta}\|fu\|_{L^\infty (\Gamma _0)}.
\end{equation}
Since $|u|\leq \delta$ in $B(\widetilde{x},r )\cap \Gamma _0$, Corollary 4.1 implies
$$
e^{-\frac{C}{r^\eta}}\leq \delta |B(\widetilde{x},r )\cap \Gamma _0|\leq \delta |\Gamma _0|
$$
or equivalently,
$$
r\leq \frac{C}{|\ln \delta |^\sigma},\; \mbox{with}\; \sigma =1/\eta .
$$
Hence, $(4.16)$ yields
\begin{equation}
\|f\|_{L^\infty (\Gamma _0)}\leq \frac{C}{|\ln \delta |^\sigma}+\frac{1}{\delta}\|fu\|_{L^\infty (\Gamma _0)},\; 0< \delta \leq \delta ^\ast .
\end{equation}
Set $\delta =e^{-s}$. Then, $(4.17)$ takes the form
\begin{equation}
\|f\|_{L^\infty (\Gamma _0)}\leq \frac{C}{s^\sigma }+e^s\|fu\|_{L^\infty (\Gamma _0)},\; s\geq s^\ast =|\ln \delta ^\ast | .
\end{equation}
We use the temporary notation $N=\|fu\|_{L^\infty (\Gamma _0)}$. The function $s\rightarrow \frac{C}{s^\sigma}+Ne^s$ attains its minimum at $\widehat{s}$ satisfying
$$
-\frac{\sigma C}{\widehat{s}^{\sigma +1}}+Ne^{\widehat{s}}=0.
$$
Using the elementary inequality $s^\varrho \leq e^{\varrho s}$, $s\geq 1$, $\varrho >0$, we obtain,
$$
\frac{C}{N}=\widehat{s}^{\sigma +1}e^{\widehat{s}}\leq e^{(\sigma +2)\widehat{s}}\;\; \mbox{if}\; \widehat{s}\geq 1.
$$
That is, 
\begin{equation}
\frac{1}{\sigma +2}\ln \frac{C}{N}\leq \widehat{s} \;\; \mbox{if}\; \widehat{s}\geq 1.
\end{equation}
But
$$
\ln \frac{C}{N}\geq \ln \frac{C}{M\|f\|_{L^\infty (\Gamma _0)}}.
$$
Therefore, there exists $\epsilon >0$ (independent on $u$ and $f$) such that $\widehat{s}\geq \max(1, s^\ast)$, provided that $\|f\|_{L^\infty (\Gamma _0)}\leq \epsilon$. When this last condition is satisfied, we can take $s=\widehat{s}$ in $(4.19)$. Taking into account $(4.18)$, $e^{\widehat{s}}=\frac{\sigma C}{N\widehat{s}^{\sigma+1}}$ and the fact that $\frac{1}{\widehat{s}^{\sigma +1}}\leq \frac{1}{\widehat{s}^\sigma}$, we easily obtain 
$$
\|f\|_{L^\infty (\Gamma _0)}\leq \frac{C}{\left|\ln \left[\|fu\|_{L^\infty (\Gamma _0)}\right]\right|^\sigma},
$$
which is the expected inequality.
\eproof

\section{Proof of  the stability theorem}
\setcounter{equation}{0}

In this section we prove Theorem 1.3. The solution of $(1.1)$ corresponding to $\lambda$ is denoted by $u(\lambda )$. Set $u^s(\lambda )=u(\lambda )-u^i$.

\medskip
We start with the following Lemma.

\begin{lemma}
Let $M>0$ be given, $\lambda \in C(\partial D)$, $0\leq \lambda \leq M$. Then there exists $R>0$, depending only on $M$ and $D$, such that $D\subset \subset B(R)$ and
\begin{equation}
|u(\lambda )|\geq 1/2,\;  |x|\geq R.
\end{equation}
\end{lemma}

{\bf Proof.} Since 
$$
|u(\lambda )(x)|=|u^i(x)+u^s(\lambda )(x)|\geq 1-|u^s(\lambda )(x)|,
$$
$(5.1)$ will follow from 
\begin{equation}
|u^s(\lambda )(x)|\leq 1/2,\;  |x|\geq R.
\end{equation}
From Green's formula of Theorem 2.4 in \cite{CK}, we have,
$$
u^s(\lambda ) (x)=\int_{\partial D}\Big[\partial _{\nu (y)} \Phi (x,y)u^s(\lambda )(y)-\partial _\nu u^s(\lambda )(y)\Phi (x,y)\Big]ds(y),\; x \in \mathbb{R}^3\setminus \overline{D}, 
$$
where,
$$
\Phi(x,y)=\frac{e^{ik|x-y|}}{4\pi |x-y|},\; x,\, y\in \mathbb{R}^3,\; x\not= y.
$$
Then, $(1.2)$ and the fact that 
$$
\partial _\nu u^s(\lambda )=-i\lambda u^s(\lambda ) -(\partial _\nu u^i+i\lambda u^i)\; \mbox{on}\; \partial D
$$ 
imply
$$
|u^s(\lambda ) (x)|\leq C\max_{y\in \partial D}\left[|\partial _{\nu (y)} \Phi (x,y)|+|\Phi (x,y)|\right],\; x \in \mathbb{R}^3\setminus \overline{D},\; |x|\geq R.
$$
A straightforward computation shows that the right hand of the last inequality tends to zero when $R$ goes to infinity. Then, $(5.2)$ follows.
\eproof

\medskip
{\bf Proof of Theorem 1.3.} Fix $R$ as in Lemma 5.1 and set $\omega =B(3R+1)$,
$$v=u(\lambda )-u(\widetilde{\lambda})= u^s(\lambda )-u^s(\widetilde{\lambda}).$$
 Recall that by estimate $(1.4)$, we have,
\begin{equation}
\|v\|_{H^3(\omega )}\leq C.
\end{equation}
Here and henceforth, $C$ is a generic constant that can depend only on $M$ and $R$. 

Let $\omega _0$ be an open subset contained in $\omega \setminus B(3R)$. Since $H^3(\omega )\subset C^{1,1/2}(\overline{\omega})$, we can apply both Propositions 1 and 2 in \cite{BD}\footnote{These two propositions are proved by similar tools to that we used in the proof of Theorem 4.1 and the main ingredient is an elliptic Carleman inequality.}. An usual argument consisting in minimizing the right hand side of estimates in Propositions 1 and 2 in \cite{BD},  with respect to the small parameter $\epsilon$, leads to the following inequality
\begin{equation}
\| v\|_{C^1(\partial D)}\leq \frac{C}{\left|\ln \left[\|v\|_{H^1(\omega _0)}\right]\right|^\kappa}\;\; \mbox{if}\; \|v\|_{H^1(\omega _0)}\leq \eta ,
\end{equation}
where, the constants $C$, $\kappa$  and $\eta$ can depend only on $M$ and $R$.
Using the interpolation inequality $\|v\|_{H^1(\omega _0)}\leq c\|v\|_{L^2(\omega _0)}^{1/2}\|v\|_{H^2(\omega _0)}^{1/2}$, $(5.3)$ and $(5.4)$, we obtain,
\begin{equation}
\| v\|_{C^1(\partial D)}\leq \frac{C}{\left|\ln \left[\|v\|_{L^2(\omega _0)}\right]\right|^\kappa}\;\; \mbox{if}\; \|v\|_{L^2(\omega _0)}\leq \eta ,
\end{equation}
 
We have by applying  Lemma 6.1.2 in \cite{Is}\footnote{This result is due to I. Bushuyev \cite{Bu}.}, where $\delta =\|u_\infty (\lambda )-u_\infty (\widetilde{\lambda})\|_{L^2(\mathbb{S}^2)}$,
$$
\| v\|_{C^1(\partial D)}\leq \frac{C}{\left|\ln \left[\delta ^{\theta (\delta )}\right]\right|^\kappa}\;\; \mbox{if}\; \delta \leq \delta _0 ,
$$
for some  constant $\delta _0 >0$. Here,
$$
\theta (\delta )=1/(1+\ln (|\ln \delta |+e)).
$$ 

Therefore, reducing $\delta _0$ if necessary,
\begin{equation}
\| v\|_{C^1(\partial D)}\leq C\left[\frac{\ln |\ln \delta |^2}{|\ln \delta |}\right]^\kappa \;\; \mbox{if}\; \delta \leq \delta _0.
\end{equation}

From the estimate in Proposition 4.1, we have
\begin{equation}
\|\lambda -\widetilde{\lambda}\|_{C(\partial D)}\leq \frac{C}{\left|\ln \left[\|(\lambda -\widetilde{\lambda})u(\lambda )\|_{C(\partial D)}\right]\right|^\sigma},
\end{equation}
if $\|\lambda -\widetilde{\lambda}\|_{C(\partial D)} \leq \epsilon$, for some $\epsilon >0$.

Or,
$$
(\lambda -\widetilde{\lambda})u(\lambda )=\widetilde{\lambda}v-\partial _\nu v.
$$
Hence,
$$
\|(\lambda -\widetilde{\lambda})u(\lambda )\|_{C(\partial D)}\leq \max (1,M)\|v\|_{C^1(\partial D)}.
$$
A combination of this last estimate, $(5.6)$ and $(5.7)$ yields
$$
\|\lambda -\widetilde{\lambda}\|_{C(\partial D)}\leq C\left|\ln \left(\frac{\ln |\ln \delta |^2}{|\ln \delta |}\right)\right|^{-\sigma} \;\; \mbox{if}\; \delta \leq \delta _0.
$$
To complete the proof, observe that the condition $\delta \leq \delta _0$ is satisfied if $\|(\lambda -\widetilde{\lambda})\|_{C (\partial D)}\leq \widetilde{\epsilon}$, for some $\widetilde{\epsilon}$, because $\lambda \rightarrow u_\infty (\lambda )$ is continuous from the set $\{ h\in C(\partial D);\; \Im h=0\; \mbox{and}\; h\geq 0\}$, endowed with the topology of $C(\partial D)$ into $L^2(\mathbb{S}^2)$ (see Appendix A).
\eproof

\appendix
\section{}
\setcounter{equation}{0}

{\bf Sketch of the proof of Theorem 1.1.} For the reader convenience, we kept the notations of \cite{CK}. Let us first recall that the fundamental solution of the Helmholtz equation $(\Delta +k^2)u=0$, with positive wave number $k$, is given as follows
$$
\Phi(x,y)=\frac{e^{ik|x-y|}}{4\pi |x-y|},\; x,\, y\in \mathbb{R}^3,\; x\not= y.
$$
We consider the single- and double-layer operators $S$ and $K$, given by
\begin{align*}
(S\varphi )(x)&=2\int_{\partial D}\Phi (x,y)\varphi (y)ds(y),\; x\in \partial D,
\\
(K\varphi )(x)&=2\int_{\partial D}\partial _{\nu (y)}\Phi (x,y)\varphi (y)ds(y),\; x\in \partial D,
\end{align*}
and the normal derivative operators $K'$ and $T$, given by
\begin{align*}
(K'\varphi )(x)&=2\int_{\partial D}\partial _{\nu (x)}\Phi (x,y)\varphi (y)ds(y),\; x\in \partial D,
\\
(K\varphi )(x)&=2\partial _{\nu (x)} \int_{\partial D}\partial _{\nu (y)}\Phi (x,y)\varphi (y)ds(y),\; x\in \partial D.
\end{align*}

As for the exterior Neumann problem, our problem is reduced to find a radiating solution $u\in C(\mathbb{R}^3\setminus D)\cap C^2(\mathbb{R}^3\setminus\overline{D})$ of the Helmholtz equation 
\begin{equation}
(\Delta +k^2)u=0\quad \mbox{in}\; \mathbb{R}^3\setminus \overline{D}
\end{equation}
satisfying the boundary condition
\begin{equation}
\partial _\nu u+i\lambda (x)u=g\quad \mbox{on}\; \partial D.
\end{equation}
Similarly to the Neumann case, we seek a solution in the form
\begin{equation}
u(x)=\int_{\partial D}\left[\Phi (x,y)\varphi (y)+i\eta \partial _{\nu (y)}\Phi (x,y)(S_0^2\varphi )(y)\right]ds(y),\; x\not \in \partial D,
\end{equation}
with a continuous density $\varphi$ and a real coupling parameter $\eta \not= 0$. The operator $S_0$ is the single-layer operator in the potential theoretic limit $k=0$. (Note that $S_0$ plays the role of a smoothing operator. We refer to \cite{CK} for more explanations.)

Next if $M_{i\lambda}$ is the multiplication operator by $i\lambda$, we will use the fact that $I+M_{i\lambda}$ is invertible. This fact is a simple consequence of the assumption that $\lambda$ is real valued.

By the results in  Theorem 3.1 in \cite{CK}, we easily prove that $u$ is a solution of $(A.3)$ provided that the density $\varphi$ is the a solution of the equation
\begin{equation}
\varphi -\left(I+M_{i\lambda} \right)^{-1}\left[K'+i\eta TS_0^2+M_{i\lambda }(S+K)\right]\varphi =-2\left(I+M_{i\lambda} \right)^{-1}g.
\end{equation}

From Theorem 3.4 in \cite{CK}, we know that the operator $(I+M_{i\lambda} )^{-1}\big[K'+i\eta TS_0^2+M_{i\lambda }(S+K)\big]$ is compact, so the Riesz-Fredholm theory is available in the space $X=C(\partial D)$. The proof will be complete if we prove that the equation $(A.4)$ with $g=0$ has only $\varphi =0$ as a solution.  

Let $u_+=u_{|\overline{D}}$ and $u_-=u_{|\mathbb{R}^3\setminus D}$. Then, $g=0$ implies that $u_-�$ is such that
$$
\partial _\nu u_-+i\lambda (x)u_-=0,\quad \mbox{on}\; \partial D.
$$
Therefore,
$$
\Im \left(\int_{\partial D}u_-\partial _\nu \overline{u}_-ds\right)=\int_{\partial D}\lambda |u|^2ds\geq 0.
$$
We deduce from Theorem 2.12 in \cite{CK} that $u_-=0$. This and the transmission conditions in Theorem 3.1 in \cite{CK} yield
$$
u_+=i\eta S_0^2\varphi ,\quad \partial _\nu u_+=\varphi \; \mbox{on}\; \partial D.
$$
Then, a simple application of Green's formula leads
$$
i\eta \int_{\partial D}|S_0^2\varphi |^2=i\eta\int_{\partial D}\varphi S_0^2\overline{\varphi}ds=\int_{\partial D}\overline{u}_-\partial _\nu u_-ds=\int_D[\nabla u|^2-k^2|u|^2]dx,
$$
whence $S_0\varphi =0$ on $\partial D$ follows. The single-layer potential $w$ with density $\varphi$ and wave number $k$ is continuous throughout $\mathbb{R}^3$  and vanishes on $\partial D$ and at infinity. Therefore, by maximum-minimum principle for harmonic functions, we have $w=0$ in $\mathbb{R}^3$ and the jump conditions in Theorem 3.1 in \cite{CK} yield $\varphi =0$.

Next, we prove estimate $(1.2)$ for the solutions of $(A.1)$-$(A.2)$. This will imply that $(1.2)$ is also valid for the solutions of $(1.1)$. To this end, we introduce the following set
$$
E=\{ M\geq 0;\; \exists\, C(M)>0\; \mbox{s.t.}\; \|u(\lambda )\|\leq C(M)\|g\|,\; \forall \, g,\, \lambda \in C(\partial D),\; 0\leq \lambda \leq M \}.
$$
Here and henceforth, $u(\lambda )$ is the solution of $(A.1)$-$(A.2)$, corresponding to $\lambda$,  and
$$
\|u(\lambda )\|=\|u(\lambda )\|_{C^\infty (\mathbb{R}^3\setminus D)},\quad \|g\|=\|g\|_{C(\partial D)}.
$$
It follows from Theorem 3.10 in \cite{CK} that $0\in E$ (corresponding to Neumann boundary condition). Let $M\in E$ and $g,\, \lambda \in C(\partial D)$, $0\leq \lambda \leq M+\epsilon$. Since
$$
\partial _\nu u(\lambda )+\frac{M\lambda}{M+\epsilon}u(\lambda )= -\frac{\epsilon \lambda}{M+\epsilon}u(\lambda ) +g,
$$
we have
$$
\|u(\lambda )\|\leq C(M)\left(\epsilon \|u(\lambda )\|+\|g\| \right).
$$
Therefore,
$$
\|u(\lambda )\|\leq \frac{C(M)}{1-C(M)\epsilon}\|g\| \;\; \mbox{if}\;\; \epsilon <1/C(M).
$$
Noting that $[0,M]\subset E$, we deduce that $E$ is an open subset of $[0,+\infty[$.

Now, let $M_k$ be a sequence in $E$, $M_k\rightarrow M$. Let $g,\, \lambda \in C(\partial D)$, $\| g\|=1$ and $0\leq \lambda \leq M$. Let $\lambda _k=M_k\lambda /M$. In view of $(4.3)$ and $(A.4)$, we easily deduce that $\lambda \rightarrow u(\lambda )$ is continuous from the set $\{ h\in C(\partial D);\; \Im h=0\; \mbox{and}\; h\geq 0\}$, endowed with the topology of $C(\partial D)$, into $C(\mathbb{R}^3\setminus D)$. Hence, there exists a positive integer $k_0$ such that $\|u(\lambda )-u(\lambda _{k_0})\|\leq 1$. Consequently,
\begin{equation}
\|u(\lambda )\|\leq \|u(\lambda )-u(\lambda _{k_0})\|+\|u(\lambda _{k_0})\|\leq 1+C(M_{k_0}).
\end{equation}
Therefore, for any $g\in C(\partial D)$, $g\neq 0$,
\begin{equation}
\|u(\lambda )\| \leq (1+C(M_{k_0}))\|g\|.
\end{equation}
We note that $(A.5)$ is trivially satisfied when $g=0$, because in this case $u(\lambda )=u(\lambda _{k_0})=0$. In conclusion, $E$ is a closed subset of $[0,+\infty[$. 

We proved that $E$ is a nonempty interval which is at the same time closed and open in $[0,+\infty [$. This implies immediately that $E=[0,+\infty[$.

With the help of Theorem 3.10 in \cite{CK}, we proceed similarly as previously to prove estimate $(1.3)$.
\eproof

\end{document}